\documentclass[11pt]{article}
\usepackage{amscd,amssymb,stmaryrd}
\usepackage{amsmath}
\usepackage{graphicx}
\usepackage{subcaption}
\usepackage[utf8]{inputenc}
\usepackage[export]{adjustbox}
\usepackage{wrapfig}
\usepackage{amstext}
\usepackage{pstricks,pst-node,pst-plot,pst-coil}
\usepackage{amsthm}
\usepackage{amsmath}
\usepackage{color}
\usepackage{mathtools}
\usepackage{hyperref}
\usepackage[english]{babel}
\usepackage{booktabs}
\usepackage{mathrsfs}
\usepackage{comment}
\usepackage{float}
\usepackage[linesnumbered,ruled,vlined, ruled,vlined]{algorithm2e}

\theoremstyle{definition}

\theoremstyle{remark}
\newtheorem*{remark}{Remark}

\title{\textcolor{black}{HEI: hybrid  explicit-implicit learning for multiscale problems}}
\author{Yalchin Efendiev\footnote{Department of Mathematics, Texas A\&M University, College Station, TX 77843, USA \& North-Eastern Federal University, Yakutsk, Russia}, ~ Wing Tat Leung\footnote{Department of Mathematics, University of California, Irvine, CA 92697, USA}, ~ Guang Lin\footnote{Department of Mathematics and Mechanical Engineering, Purdue University, West Lafayette, IN 47906, USA}, ~Zecheng Zhang \footnote{Department of Mathematics, Purdue University, West Lafayette 47906, IN, USA}}

\begin{document}
\maketitle

\begin{abstract}
Splitting method is a powerful method to handle application problems by splitting physics, scales, domain, and so on. Many splitting
algorithms have been designed for efficient temporal discretization. 
In this paper, our goal is to use temporal splitting concepts in
designing machine learning algorithms and, at the same time,
help splitting algorithms by incorporating data and speeding them up. \textcolor{black}{Since the spitting solution usually has an explicit and implicit part, we will call our method hybrid explicit-implict (HEI) learning.}
We will consider a recently introduced multiscale splitting algorithms, where the multiscale problem is solved on a coarse grid. To approximate the dynamics, only a few degrees of freedom are solved implicitly, while others explicitly. 
This splitting concept allows identifying degrees of freedom that need implicit treatment. In this paper, we use this splitting concept in machine learning and propose several strategies. First, the implicit part of the solution
can be learned as it is more difficult to solve, while the explicit part can be computed. This provides a speed-up and data incorporation for splitting approaches. Secondly, one can design a hybrid neural network architecture because handling explicit parts requires much fewer communications among neurons and can be done efficiently. Thirdly, one can solve the coarse grid component via PDEs or other approximation methods and construct \textcolor{black}{simpler neural networks for the explicit part of the solutions}. \textcolor{black}{We discuss these options
and implement one of them by interpreting it as a machine translation task.
This interpretation of the splitting scheme successfully enables us using the Transformer 
since it can perform model reduction for multiple time series and learn the connection between them.
We also find that the splitting scheme is a great platform to predict the coarse solution with insufficient information of the target model: the target problem is partially given and we need to solve it through a known problem which approximates the target. Our machine learning model can incorporate and encode the given information from two different problems and then solve the target problems.}
We conduct four numerical examples and the results show that our method is stable and accurate.

\end{abstract}

\section{Introduction}

Splitting methods are used for many applications, where splitting physics, scales, domains, and so on, are considered. These approaches allow simplifying
problems and solving them separately and then coupling. Recently, novel splitting approaches are proposed for multiscale problems. The main idea of these approaches is to identify the coarse grid components of the solution and treat them
implicitly, while treating other components explicitly. In particular, 
the stability analysis shows that one needs to choose spaces appropriately
such that the first space contains all components associated with high contrast and multiscale features of the solution. Using this splitting, we show that 
one can choose the time step that scales as the coarse mesh size
and the time step is independent of the contrast.

This approach, on the one hand, allows simplifying the temporal discretization
and using explicit schemes with some minimal time step constraints. On the 
other hand, these approaches allow designing neural network architectures
and/or using them to speed up the solutions. Before discussing it,
we briefly describe some main ideas of splitting algorithms.

We assume that the forward problem is given by 
\[
u_t = f(u),
\]
where $f(u)$ is a multiscale differential operator. 
To achieve a contrast-independent stability, 
 splitting $u=u_1+u_2$ is performed, where $u_1$ is a coarse-grid 
approximation and uses much fewer degrees of freedom, while $u_2$ is
other degrees of freedom. To have a stability, one needs certain
conditions on the space for $u_2$, such as $\|u\|_{L^2}/\|u\|_a$ is bounded
below independent of the contrast, where $\|\cdot\|_a$ is a norm associated
with the operator $f$ (roughly speaking). To achieve this, one needs
to include all high-contrast multiscale features in $u_1$. This construction
is described in \cite{chung2021contrast} 
and uses GMsFEM and NLMC concepts
developed in \cite{chung2018constraintmixed,NLMC}. 
Other multiscale methods have been developed
for various applications \cite{eh09,hw97,jennylt03,chung2016adaptiveJCP,MixedGMsFEM,WaveGMsFEM,chung2018fast,GMsFEM13,chung2018constraint, chung2018constraintmixed,NLMC,ee03,henning2012localized,rk07,skr06, chetverushkin2021computational,chung2021computational}, 
which we do not mention. \textcolor{black}{In this paper, we will focus on 
NLMC and GMsFEM.}

\textcolor{black}{One of the contributions} of this paper is the use of machine learning to 
accelerate the computations and discuss a possible design for machine
learning \cite{chung2021multi, zhang2020learning, chung2020multi, liu2021deep, leung2021nh}. We again comment that the computation of
$u_1$ is more computationally difficult compared to $u_2$ as it uses
implicit discretization.
For this reason, we propose the following strategies.

\begin{itemize}

\item 1) We learn and predict $u_1$ by designing appropriate 
machine learning techniques, while we compute $u_2$ based on PDEs. 

\item 2) We design a hybrid neural network, one for $u_1$ and the other
for $u_2$, and take advantage of the fact that $u_2$ is explicit.

\item 3) We compute $u_1$ based on data, reduced-order models, ,.. and
learn $u_2$ due to its sparse forward map structure.

\end{itemize}

\textcolor{black}{All strategies consist of learning either $u_1$ implicitly or $u_2$ explicitly, we hence call our methodology: hybrid explicit-implicit learning (HEI).}
Next, we briefly discuss these cases. For the first case, our motivation
stems from the fact that the computation of $u_1$ is expensive as it 
is implicit and nonlinear (in general). Thus, for its computation, 
machine learning techniques can be used to predict $u_1$, while we compute
$u_2$ since its computation is explicit.
\textcolor{black}{More precisely, we can train a machine learning model based on a small portion
of the entire solution in time and then alternatively predict $u_1$ and compute $u_2$.}

For the second case, our main idea is to design a hybrid neural 
network.
The computation for $u_2$ uses sparse information and needs 
fewer parameters for identification.
\textcolor{black}{However, if the problem is nonlinear, 
one needs to construct the multiscale basis space for $u_1$ and $u_2$ at each iteration.
This increases the computational cost in computing $u_2$ by the splitting scheme,
hence training and predicting  $u_2$
is a great option.}

For the third case, 
our motivation comes from the fact that $u_1$ can be observed
or extracted from the solution dynamics (e.g., approaches like DMD).
\textcolor{black}{On the other hand, the evolution of $u_2$ uses sparse information
and it may not be difficult to train a machine learning model to predict $u_2$.} 
We then can use approximate dynamics (learned or a priori formulated) to compute
$u_2$.

In this paper, we mostly focus on the first approach. In particular,
we use some special neural network architectures for it, which are based
on the Transformer.
Transformer is a machine learning structure which is originally designed to handle language tasks \cite{vaswani2017attention}.
It can be regarded as a machine learning model which can encode the long time dependency of the time series \cite{dai2019transformer, vaswani2017attention, velivckovic2017graph, zhang2019self, chung2021multi}.
One hence can apply it to perform model reduction in time \cite{chung2021multi}. 
Our problem is to predict the dynamics of $u_1$ or $u_2$ which 
are of high dimension, so it is natural to use Transformer to learn the time dependency meanwhile performing the model reduction.
Furthermore, Transformer can be used in handling the machine translation task, i.e., 
it can deal with the problems which have inputs of different space (languages).
This motivates us using Transformer to learn $u_1$ because the dynamics of $u_1$ 
depends on both $u_1$ and $u_2$ which are in different spaces.
We conduct four experiments by using Transformer as the key of the machine learning model
and the results show that the prediction \textcolor{black}{is stable and is accurate}.

Another contribution of this work is the \textcolor{black}{coarse solution prediction with incomplete information,
that is, we find that the splitting scheme is a good platform to find the coarse solution of a model, 
but only the partial information of the target model is given.
For simplicity, we will call this `coarse solution assimilation'
in the remaining of this work. }
There are many real life problems whose
solution $v$ (different from $u$)
can be decomposed as $v_{1}+v_{2}$, where $v_1$ can capture the main features of the model
while $v_2$ is the correction.
In the real scenario, only incomplete information of $v_1$ is given and information about $v_2$ is missing.
If one wants to recover $v$, one way is to solve a model which approximates the target model, 
i.e., we know the solution $\Tilde{v} = \Tilde{v}_1+\Tilde{v}_2$ of a model which is closed to the target.
The problems like this can be referred to the coarse solution problem.
\textcolor{black}{This problem is not easy to solve,} but 
\textcolor{black}{we realize that the splitting scheme can be formulated as a coarse solution assimilation process} and 
we can use our proposed machine learning method to solve the problem easily.

As we have discussed before, $u_1$ is a coarse-grid approximation and can capture all 
high-contrast multiscale features.
This can be regarded as $v_1$ in the coarse solution assimilation.
Moreover, $u_2$ uses other degree of freedom and is the correction $v_2$.
%We hence believe that the the data assimilation is similar to the splitting scheme.
Now the question is how to solve the problem, i.e., find all $u_1+u_2$ 
if it is only given an incomplete $u_1$ and $\Tilde{u}_2$ which comes from 
an approximating model.

We find that our machine learning model is able to handle this problem.
We can train a model which predicts the future $u_1$ with a small portion of 
given $u_1$ and $\Tilde{u}_2$.
The difficulty is to handle inputs which originate from different spaces ($u_1$ and $\Tilde{u}_2$) 
and also learn the time dependency. This includes encoding the $u_1$, $\Tilde{u}_2$ and combining
 $u_1$ and $\Tilde{u}_2$.
We find that our proposed method can handle this problem and give a reliable prediction.
It should be noted that the strategies mentioned previously can all be incorporated into 
the coarse solution assimilation, \textcolor{black}{that is, we can propose hybrid explicit-implicit (HEI) learning strategies to solve the problem.}
For the demonstration purposes, we will only consider the first case.

In the paper, \textcolor{black}{we present several numerical results basing on our hybrid explicit-implicit (HEI) learning framework.}
We solve both linear and nonlinear examples and their applications in coarse solution assimilation. 
We observe that the solutions are stable and both $L_2$ and energy relative errors are small.

The rest of the paper is organized as follows.
In Section \ref{splitting_sec}, we briefly review the splitting scheme, the motivation of the machine learning and the basic machine learning model are introduced in Section \ref{machine_learning_sec}. We present our proposed method, \textcolor{black}{hybrid explicit-implicit (HEI) learning}, in Section \ref{proposed_sec} and introduce the splitting-based coarse solution assimilation in Section \ref{data_assimilation_sec}. Finally, in Section \ref{numerical_sec}, we verify our proposed method by conducting 4 numerical experiments. The results are also shown in this section.

\section{Preliminaries}
\label{splitting_sec}

\textcolor{black}{In this section, we will review the partially explicit splitting scheme \cite{chung2021contrast, Chung2021ContrastindependentPE} which is designed to tackle the multiscale high contrast problems. The scheme is based on the multiscale finite element methods \cite{chung2018constraintmixed,NLMC} and the stability of the scheme is independent of the contrast ratio of the problem.}
We consider the following equation
\begin{align}
u_{t}=-f(u). \label{eq1}
\end{align}
In general, we can assume that
$f=\cfrac{\delta F}{\delta u}$, \textcolor{black}{the variational derivative of energies $F(u) = \int E_1(u)$}
, and is contrast dependent and nonlinear (or linear) (i.e.,
$f$ introduces stiffness in the system).
We refer to \cite{Chung2021ContrastindependentPE} 
for more details about assumptions on
$F$ for well-posedness, weak formulation, and stability.

To solve the problem, a standard method is finite element
approach. We can consider the numerical solution $u_{H}(t,\cdot)\in V_{H}$
satisfies 
\begin{align}
\textcolor{black}{(u_{H,t},v)=-(f(u),v)}\;\forall v\in V_{H}, \label{eq:CEM_problem}
\end{align}
where $V_{H}$ is a finite element space and \textcolor{black}{$(\cdot, \cdot)$ denotes the inner product in $V_H$.}

\textcolor{black}{We consider the solution splitting and split $V_{H}$ as
a direct sum of two subspace $V_{H,1}$ and $V_{H,2}$, namely, 
$V_{H}=V_{H,1}\oplus V_{H,2}$. In practice, $V_{H, 1}$ and $V_{H, 2}$ are the coarse multiscale finite element spaces.}
The finite element solution is then satisfying
\begin{align*}
(u_{H,1,t}+u_{H,2,t},v_{1})+(f(u_{H,1}+u_{H,2}),v_1) & =0\;\forall v_{1}\in V_{H,1},\\
(u_{H,1,t}+u_{H,2,t},v_{2})+(f(u_{H,1}+u_{H,2}),v_2) & =0\;\forall v_{2}\in V_{H,2},
\end{align*}
where $u_{H}=u_{H,1}+u_{H,2}$. 
\textcolor{black}{For the time discretization, we  consider the temporal splitting 
introduced in \cite{chung2021contrast} for linear problems
and in \cite{Chung2021ContrastindependentPE} 
for nonlinear problems. 
That is, we can use a partially explicit time
discretization.} For example, we can consider 
\begin{align*}
& (\cfrac{u_{H,1}^{n+1}-u_{H,1}^{n}}{\tau}+\cfrac{u_{H,2}^{n}-u_{H,2}^{n-1}}{\tau},v_{1})+(f(u_{H,1}^{n+1}+u_{H,2}^{n}),v_{1})  =0\;\forall v_{1}\in V_{H,1},\\
& (\cfrac{u_{H,1}^{n}-u_{H,1}^{n-1}}{\tau }+\cfrac{u_{H,2}^{n+1}-u_{H,2}^{n}}{\tau },v_{2})+\textcolor{black}
{(f(\omega u_{H,1}^{n+1}+ (1-\omega)u_{H,1}^{n+1}+u_{H,2}^{n}),v_{2})}  =0\;\forall v_{2}\in V_{H,2},
\label{splitting_nl_splitting_scheme}
\end{align*}
\textcolor{black}{where $\tau$ is the time step and $\omega\in[0, 1]$ is a constant. For the linear case, the splitting scheme is}
\begin{align}
    & (u_{H, 1}^{n+1}, v) = (u_{H, 1}^{n}, v) - (u_{H, 2}^{n}-u_{H, 2}^{n-1}, v)-\tau a(u_{H, 1}^{n+1}+u_{H, 2}^{n}, v), v\in V_{1, H},\\
    & (u_{H, 2}^{n+1}, v) = (u_{H, 2}^{n}, v) - (u_{H, 1}^{n}-u_{H, 1}^{n-1}, v)-
    \tau a\big( (1-\omega)u_{H, 1}^{n} +\omega u_{H, 1}^{n+1}+u_{H, 2}^{n}, v \big), v\in V_{2, H},
    \label{splitting_splitting_scheme}
\end{align}
\textcolor{black}{where $a(\cdot, \cdot)$ is the bilinear form associated with $f$.
One example for $f$ is the parabolic problem, i.e., $f(u) = \nabla\cdot(\kappa(x)\nabla u)$, where $\kappa$ is the contrast permeability field. The bilinear form is then given by $a(u, v) = \int_{\Omega}\kappa(x) \nabla u\nabla vdx$ for $u, v\in V_H$ and $\Omega$ is the domain of the problem. For the stability and convergence concerns of the scheme, we refer to \cite{chung2021contrast, Chung2021ContrastindependentPE}.
}
% In many applications, coarse-grid components can be predicted \textcolor{red}{when} there are data available. For example, in data assimilation (e.g., \cite{azouani2014continuous,bessaih2015continuous}), \textcolor{red}{one of the models induces the splitting scheme as follow:}
% \begin{align}
%     & (u_{H, 1}^{n+1}, v) = (u_{H, 1}^{n}, v) - (u_{H, 2}^{n}-u_{H, 2}^{n-1}, v)-\tau a(u_{H, 1}^{n+1}+u_{H, 2}^{n}, v) + \beta (u_{H,1}^{obs,n+1}-u_{H, 1}^{n+1}, v), v\in V_{1, H},\\
%     & (u_{H, 2}^{n+1}, v) = (u_{H, 2}^{n}, v) - (u_{H, 1}^{n}-u_{H, 1}^{n-1}, v)-
%     \tau a\big( (1-\omega)u_{H, 1}^{n} +\omega u_{H, 1}^{n+1}+u_{H, 2}^{n}, v \big), v\in V_{2, H},
%     \label{splitting_scheme}
% \end{align}
% for some $\beta$ which imposes the constraints given by observation data. 
% From the splitting algorithm, we see that this part of the solution
% is difficult to compute and we propose to use machine learning for
% prediction of ``coarse-grid'' components and 
% ``residual model'' for predicting possibly important fluctuations/corrections.

% \begin{figure}[H]
% \centering
% \includegraphics[scale = 0.4]{graphs/chart1.png}
% \caption{}  
% \label{fig1}
% \end{figure}

\section{Machine learning}
\label{machine_learning_sec}
\subsection{Motivation}
\label{sec_ml_motivation}
One of the great benefits of the splitting scheme is the contrast independent property.
However, it doubles the computations since we need to iterate twice in time.
Machine learning is widely used in solving mathematics problems. 
Compared to the implementation of splitting numerically, the deep learning approach may be more computational friendly.
In particular, researchers need to solve a matrix inverse problem in the first equation of the splitting scheme (\ref{splitting_splitting_scheme}).
However, one only needs to evaluate a sequence of matrix multiplications in the deep learning framework.
\textcolor{black}{Motivated by this, we propose the hybrid explicit-implicit (HEI) learning method. }

\textcolor{black}{The basic idea is to first compute $\{u_{H, 1}^n, u_{H, 2}^n\}_{n= 1}^M $, where $M<<N$ and $N$ is the total number of time steps in time, then we can train machine learning models to predict either 
$u_{H, 1}^{n}$ or $u_{H, 2}^{n}$, where $n = M+1, ..., N$.
We hence propose three strategies.
}
\begin{enumerate}
    \item Train a model to predict $u_{H, 1}^n$ while computing the $u_{H, 2}^n$ according to the splitting scheme.
    The motivation of this option is that $u_{H, 1}^n$ is implicit and is not easy to compute.
    \item Train a model to predict $u_{H, 2}^n$ while computing $u_{H, 1}^n$ by some other methods. The motivation of this idea is that the evolution of $u_{H, 2}^n$ is explicit and should be easier to learn.
    \item Train two models to predict both $u_{H, 1}^n$ and $u_{H, 2}^n$. Learning both solutions is also beneficial in particular when one needs to solve a nonlinear problem.
\end{enumerate}

We have discussed the benefits of \textcolor{black}{three hybrid explicit-implicit (HEI) learning strategies} before in the Introduction and in this work, 
we will focus and elaborate on the first strategy, i.e., 
we will introduce the machine learning model and discuss the coarse solution assimilation idea based on this strategy in section (\ref{data_assimilation_sec}). The other two strategies also rely on machine learning models 
and our proposed model, training and testing methods can be applied directly.

Now, let us dive into the first model problem.
More precisely, we are going to train a model which predicts $u_{H, 1}^{n+1}$ using a small portion of the computation data in $[0, t]$ where $t<<T$,
then in the remaining greater computational domain $[t, T]$, we will apply the trained model to predict $u_{H, 1}^{n+1}$,
this value is then used in the second equation to obtain $u_{H, 2}^{n+1}$.
It should be noted that 
the proposed scheme is built based on the splitting scheme (\ref{splitting_splitting_scheme}),
our model mimics the first equation and predicts $u_{H, 1}^{n+1}$.
For $u_{H, 2}^{n+1}$ evaluation, we still follow the second equation since it is also a forward process and does not involve intense computation.
The proposed method can be summarized as follows.
\begin{enumerate}
    \item Compute $u_{H, 1}^{n}$ and $u_{H, 2}^{n}$ up to $n = M$, where $M<<N$ and $N$ is the total simulation step.
    The collection $\{u_{H, 1}^{n}, u_{H, 2}^{n}\}_{n = 1}^M$ will be denoted as $\mathcal{D}$.
    \item Train a model model $\mathcal{N}$ which predicts $u_{H, 1}^{n+1}$ using the training dataset $\mathcal{D}$. More precisely, $\mathcal{N}$ is a map from
    $\big\{ \{u_{H, 2}^{n-i+1}\}_{i = 1}^{n_e}, \{u_{H, 1}^{n+i-1}\}_{i = 1}^{n_d} \big\}$ to $u_{H, 1}^{n+1}$, where $n_e$ and $n_d$ are the input dimension of the encoder and decoder. Theoretically, $n_e = 2$ and $n_d = 2$, we will detail this in the next section.
    \item From $M$ to $N$, predict $u_{H, 1}^{n+1}$ by the model $\mathcal{N}$ and compute $u_{H, 2}^{n+1}$ by the equation iteratively until the terminal time.
\end{enumerate}

\subsection{Self attention and Transformer}
We find that the splitting scheme is similar to the machine language translation task in nature.
We will apply Transformer to design the machine learning model $N$.
In this section, let us first review the attention and Transformer.

\textbf{Self attention.} 
Self attention \cite{vaswani2017attention, chung2021multi} is a model to handle time series data.
The input is a time series consisting of multiple time steps.
The module can be regarded as a model reduction \cite{chung2021multi} and will output
a new series with a reduced dimension.
It encodes a long context into a smaller size representation \cite{dai2019transformer}.
\textcolor{black}{In particular, this encoding is able to relate different positions
of this sequence \cite{vaswani2017attention} and concentrate on the input's most relevant parts to perform the task \cite{velivckovic2017graph}.
The response at a time step position is the weighted sum of the features (feature: reduced order model for a single time step) at
all positions
the long range, multi-level dependencies
can be enforced \cite{zhang2019self}.}
Self-attention has proven to be useful for tasks such as machine
reading and learning sentence representations \cite{velivckovic2017graph}.

\begin{figure}[H]
\centering
\includegraphics[scale = 0.4]{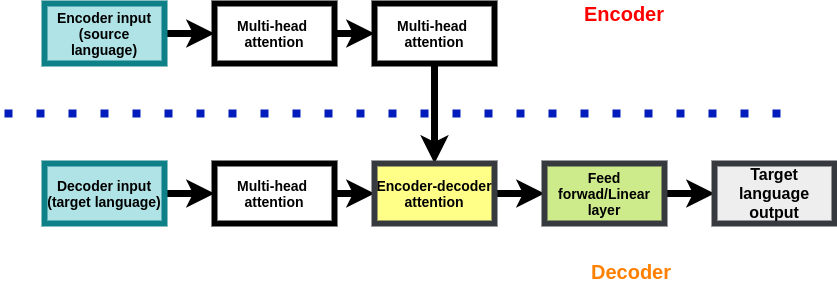}
\caption{Transformer. The top of the image is the encoder and the bottom is the decoder. Multi-head attention (white) is the stack of self-attention blocks. The encoder-decoder attention (yellow) helps the decoder focus on appropriate places in the source input, it is the key to learn the connections between the encoder and decoder. We skip the positional encoding, softmax and other blocks which are not the focus of this work. Please refer to \cite{vaswani2017attention} for details.}  
\label{ml_transformer}
\end{figure}

\textbf{Transformer.} Transformer is an architecture for transforming one sequence into another one (eg. from English to Germany).
It is built on the self-attention module and likes the most transduction models, it consists of an encoder and a decoder which are the stack of self attention modules. Please check Figure (\ref{ml_transformer}) and refer to \cite{vaswani2017attention} for the details of Transformer.

Transformer is able to build connections between the encoder and decoder input and output \cite{vaswani2017attention}. More precisely, the encoder output is connected to the decoder by an encoder-decoder attention block, this block is able to learn the source and target language connections.

We review the inference phase of Transformer. 
This is an iterative process which is similar to the time evolution of the numerical scheme: the prediction of the next $u_{H, 1}^n$ relies on the previous predictions.
The testing has the following steps.
\begin{enumerate}
    \item Input the full encoder sequence (source language, e.g., English). However,  
    the decoder input is an an empty sequence with a beginning token.
    The output is a sequence. The first element in the output sequence will be the first word of the target sentence (e.g., Germany). 
    This word will also be used in the prediction of the following words in the target sentence.
    \item The first word will be inserted into second place of the decoder input sequence. Together with the entire encoder input (the whole English sentence),  Transformer will output
    a new sequence whose the first element will be the second word of the target sentence.
    \item Repeat this process until the end token is predicted and one then gets the whole translated sentence.
\end{enumerate}
The translation process is also iterative, i.e., the next word relies on the previous predictions.
However, Transformer is capable to handle the cumulative error or correct the mistakes made in the early stage of the prediction.
This motivates us to use Transformer to model the splitting scheme.

\section{\textcolor{black}{Hybrid explicit-implicit (HEI) learning}}
\label{proposed_sec}
In Section (\ref{sec_ml_motivation}), \textcolor{black}{we have introduced the implementation of the hybrid explicit-implicit (HEI) learning strategies}, however, we do not give the details about how to incorporate the splitting with Transformer. In this section, we are going to introduce the motivation of choosing Transformer and give the structure of the network.

\subsection{Motivation of using Transformer}
It is very natural to design a machine learning network based on the splitting scheme directly.
As we have discussed before, the computation involved in the second equation is only the matrices multiplications. We hence focus on the first equation and theoretically, we need to approximate the matrix inversion and multiplication by the network. 
There should be some model design principles and 
we summarize the principles based on the first equation as follows.

\begin{enumerate}
    \item \textbf{Principle I}. To get $u_{H, 1}^{n+1}$, one needs $u_{H, 1}^{n+1}$ and $\{u_{H, 2}^{n}, u_{H, 2}^{n-1} \}$.
They come from different space $V_{H, 1}$ and $V_{H, 2}$, the model should be able to encode two different time series independently and learn the interaction between the corresponding reduced order models.

\item \textbf{Principle II}. The input from $V_{H, 2}$ has two time steps $\{u_{H, 2}^{n}, u_{H, 2}^{n-1} \}$, the model should be capable of capture the dependency of the multiple-step input.

\item \textbf{Principle III}. Due to the stability concern, the accumulative error should be controlled in a proper way, that is, if we predict $u_{H, 1}^{n+1}$ by the model and compute $u_{H, 2}^{n+1}$ by the second equation iteratively, the solutions cannot blow up.

\end{enumerate}

Recall the self-attention and Transformer, they satisfy all design principles. 
$\{u_{H, 1}^{n-i+1}\}_{i = 1}^{n_e} $ can be regarded as the target language and input of the decoder, and $ \{u_{H, 2}^{n+i-1}\}_{i = 1}^{n_d}$ is the source language and input of the encoder. We then can summarize how Transformer handles the design difficulties as follows.
\begin{enumerate}
\item \textbf{Principle I.} According to \cite{vaswani2017attention}, Transformer is able to build the global dependency between the input and output.
This is achieved by the encoder-decoder attention 
which learns the influence of the source language ($\{u_{H, 2}^{n-i+1}\}_{i = 1}^{n_e}$) on the target language ($\{u_{H, 1}^{n+i-1}\}_{i = 1}^{n_d}$).

\item \textbf{Principle II.} The self attention mechanism built in Transformer is able to learn the dependency and influence among different positions (time steps).

\item \textbf{Principle III.} As reviewed before, Transformer is auto-regressive, i.e., it uses the last predicted word ($u_{H, 1}^{n+1}$) in the target sentence as the additional input to infer next word  ($u_{H, 1}^{n+2}$) in the sentence. Transformer has been successful in various language tasks, we hence believe it can deal with the accumulative errors.
\end{enumerate}

\subsection{Proposed structure and implementation}
We give the details of the network training and testing method in this section. The splitting scheme is similar to the language translation task, we hence use the similar training and testing strategy as the language task.

\begin{figure}[H]
\centering
\includegraphics[scale = 0.35]{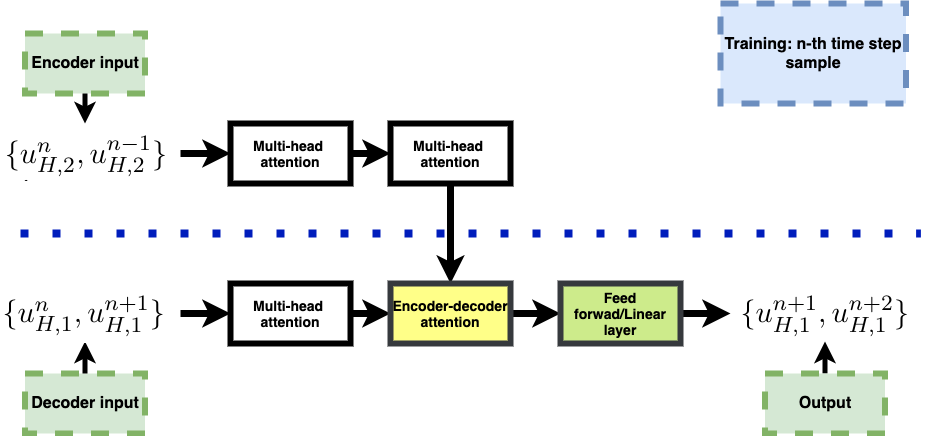}
\caption{Splitting Transformer model training phase. We show one training sample at time step $n$. Note that we set $n_d = n_e = 2$, so the model is mimicking the first equation of the splitting scheme.
}   
\label{pp_train}
\end{figure}

\begin{figure}[H]
\centering

\includegraphics[scale = 0.35]{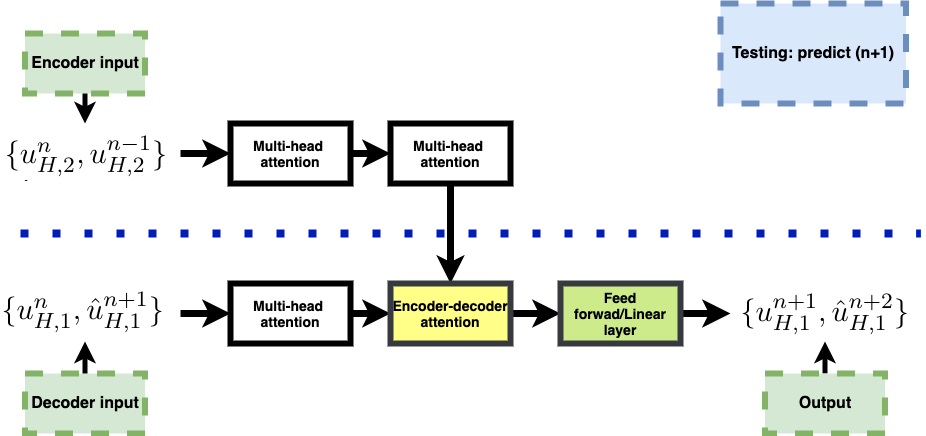}
\caption{Splitting Transformer model testing phase \textcolor{black}{($n_e = n_d = 2$)}. We show the prediction of $u_{H, 1}^{n+1}$, which will be used in calculating 
$u_{H, 2}^{n+1}$.
The input $\hat{u}_{H, 1}^{n+1}$ comes from the second entry of the last step prediction,
$\hat{u}_{H, 1}^{n+2}$ will be used as the input of the next prediction of $u_{H, 1}^{n+2}$.
We design the training testing flow in this way to mimic the implicit scheme of the first equation.
}   
\label{pp_test}
\end{figure}

\textbf{Structure.}
The encoder input is $\{u_{H, 2}^{n-i+1}\}_{i = 1}^{n_e}$,
this can be regarded as the source language with length $n_e$.
The decoder input is the $\{u_{H, 1}^{n+i-1}\}_{i = 1}^{n_d}$,
we take this as the target language. 
The structure of the network is shown in Figure (\ref{pp_train}) and Figure (\ref{pp_test}).

\textbf{Training.}
During the training,  $\{u_{H, 2}^{n-i+1}\}_{i = 1}^{n_e}$ and
$\{u_{H, 1}^{n+i-1}\}_{i = 1}^{n_d} $ will be inputs, the output of the network will be $\{u_{H, 1}^{n+i}\}_{i = 1}^{n_d}$.
The training of the network is demonstrated in Figure (\ref{pp_train}).

\textbf{Testing.}
This is an iterative process and at the step $n$, the prediction $u_{H, 1}^{n+1}$ of the network will be used in the calculating $u_{H, 2}^{n+1}$ in the second equation. We summarize the process as follows and demonstrate the prediction of $u_{H, 1}^{n+1}$ in Figure (\ref{pp_test}).
\begin{enumerate}
    \item \textcolor{black}{At step $n$, the encoder input is $\{u_{H, 2}^{n-i+1}\}_{i = 1}^{n_e}$ and the decoder input is$\{u_{H, 1}^{n}, \hat{u}_{H, 1}^{n+1}, ..., \hat{u}_{H, 1}^{n+n_d-1}\}$;}
    the output of the network is $\{u_{H, 1}^{n+1}, \hat{u}_{H, 1}^{n+2}, ..., \hat{u}_{H, 1}^{n+n_d} \}$.
    \item The first entry $u_{H, 1}^{n+1}$ will be used in the second equation to calculate $u_{H, 2}^{n+1}$;
    \textcolor{black}{while $\hat{u}_{H, 1}^{n+2}, ..., \hat{u}_{H, 1}^{n+n_d}$ will be used in predicting $u_{H, 1}^{n+2}$ and future $u_{H, 1}$ sequence.}
    \item Set $n = n+1$ and repeat step 1 and step 2 until reaching the terminal time.
\end{enumerate}
If we follow the splitting scheme,  $n_e = 2$ and $n_d = 2$.

\begin{remark}
It should be noted that the first equation is implicit in $u_{H, 1}^{n+1}$ at step $n$ and our scheme is mimicking the process. At step $n+1$,
the target of the prediction is $u_{H, 1}^{n+2}$, however, 
this term has already been predicted in the last iteration $n$
and is used in the input of the decoder at time $n+1$.
\textcolor{black}{If $n_d = 2$, denoted as $\hat{u}_{H, 1}^{n+2}$, this term is the second entry of the output at 
step $n$.}
We do not use $\hat{u}_{H, 1}^{n+2}$ in the calculating of
$u_{H, 2}^{n+2}$, instead we use \textcolor{black}{$u_{H, 1}^{n+2}$ which is the first entry of the output at step $n+1$}. Please check Figure (\ref{pp_test}) for the illustration.
\end{remark}

\subsection{Coarse solution assimilation}
\label{data_assimilation_sec}
\textcolor{black}{In this section, we will discuss applications of our method. 
We find that the splitting scheme can 
be used to solve the target problems, however, a part of the
 model is implicitly given. 
For simplicity, we will call this part of the 
 problem the coarse solution assimilation.}
There are many real life applications whose target model is not given or partially given.
Researchers then need to solve the target model based on a known model which approximates the target model. 
This is an interesting problem and we will explain this model and our method in this section.

We assume that the target solution $u_{\kappa}$ of the target equation with permeability $\kappa(x)$ can be decomposed as 
$u_{\kappa}^n = u_1^n+u_2^n, n = 1,..., N$, where $u_1^n$ can be regarded as the coarse model while $u_2^n$ is the correction, $n$ here indicates the time step.
We are lack of $\kappa(x)$ information but still want to evaluate $u_{\kappa}^n$, where $n = 1,..., N$.
We claim that if we know $\Tilde{\kappa}(x)$ which approximates $\kappa(x)$,  and a small piece of $\{u_{1}^n\}_{n = 1}^M$, where $M<<N$, then $u_{\kappa}^n$ can be recovered by our method. Intuitively, we propose to learn a model
that predicts the coarse solution $u_1^n$ of the target equation.
Moreover, the coarse solution is automatically
corrected by $\Tilde{u}_2^n$ from a known equation.

\textcolor{black}{
Furthermore, in the framework of the splitting scheme, we have the following discretization for the linear case
\begin{align}
    & (u_{H, 1}^{n+1}, v) = (u_{H, 1}^{n}, v) - (\Tilde{u}_{H, 2}^{n}-\Tilde{u}_{H, 2}^{n-1}, v)-\tau a(u_{H, 1}^{n+1}+\Tilde{u}_{H, 2}^{n}, v) , v\in V_{1, H},\\
    & (\Tilde{u}_{H, 2}^{n+1}, v) = (\Tilde{u}_{H, 2}^{n}, v) - (u_{H, 1}^{n}-u_{H, 1}^{n-1}, v)-
    \tau a\big( (1-\omega)u_{H, 1}^{n} +\omega u_{H, 1}^{n+1}+\Tilde{u}_{H, 2}^{n}, v \big), v\in V_{2, H},
    \label{splitting_assimilation_2}
\end{align}
where $\Tilde{u}_{H, 2}$ is the given observation which is the contrast independent component of an approximating solution,
this solution may be generated from a permeability $\Tilde{\kappa}$ which is close to the target permeability field $\kappa$ (not given).
In this case, we can still predict the coarse grid component $u_{H, 1}^{n}$ and the given $\Tilde{u}_{H, 2}^{n}$ can be regarded as the correction of the coarse component solution. Since $\Tilde{u}_{H, 2}^n$ is not the real correction of $u_{H, 1}^n$, the scheme is not accurate.
However, the machine learning can capture the relation between $\Tilde{u}_{H, 2}^n$ and $u_{H, 1}^n$. More precisely, 
the benefit and the target of machine learning is that
 the model is supposed to learn and predict all $u_{H, 1}^{n}$ with a small piece of $u_{H, 1}^{n}$ and $\Tilde{u}_{H, 2}^{n}$, such that $\Tilde{u}_{H, 2}^{n}$ is able to improve 
$u_{H, 1}^{n}$ and the combined solution $u_H^n = u_{H, 1}^n$+$\Tilde{u}_{H, 2}^n$ is accurate and stable. In this work, we are going to study this example.
}
\begin{figure}[H]
\centering
\includegraphics[scale = 0.4]{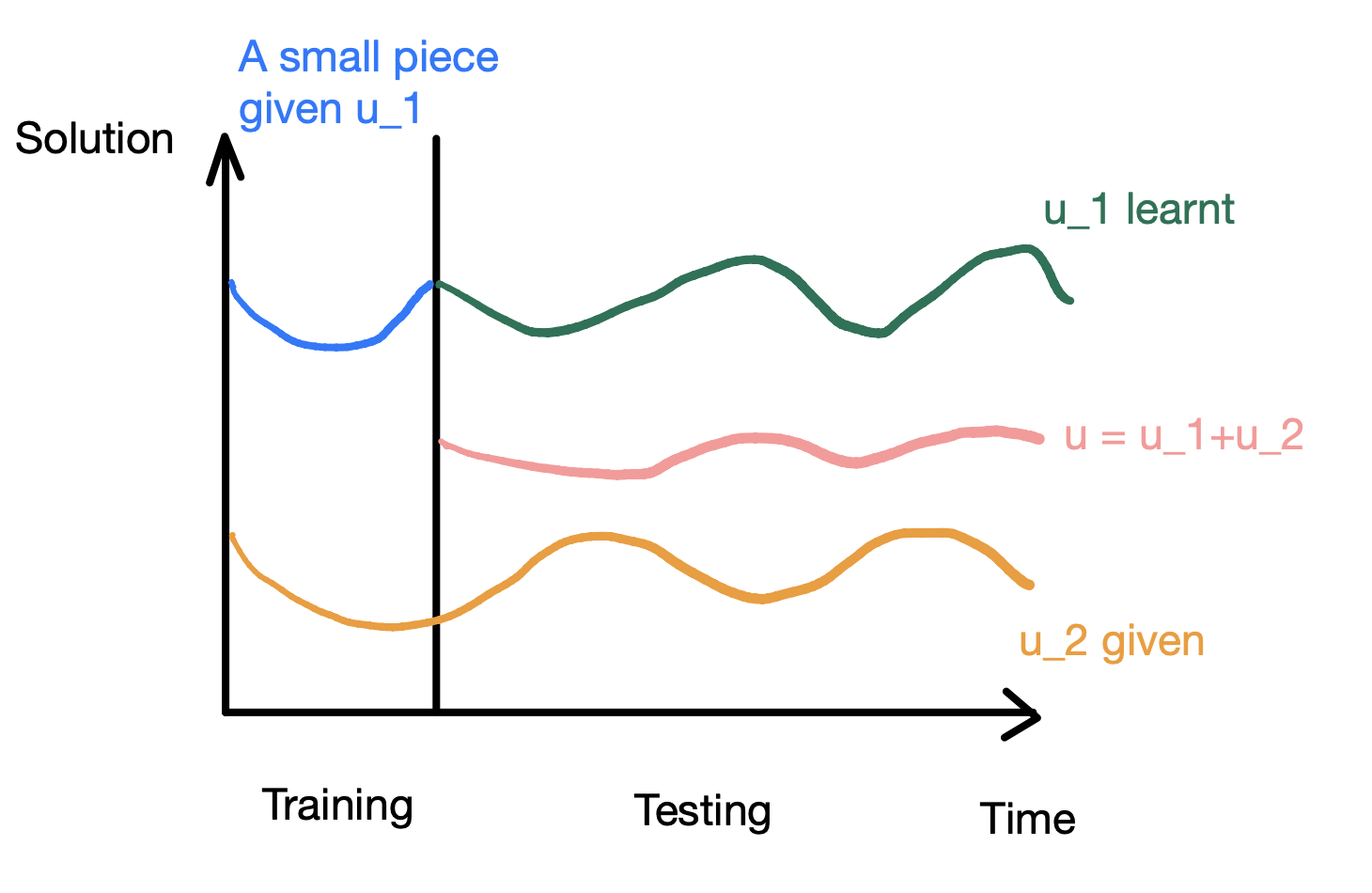}
\caption{\textcolor{black}{Demonstration of the coarse solution assimilation. We assume a small piece of $u_1$ and the entire $u_2$ are given. $u_2$ may come from a model which approximates the target model. $u = u_1+u_2$ is the combined solution, in which the learnt $u_1$ captures the multiscale properties and $u_2$ is the correction. The machine learning model should be able to learn $u_1$ basing on $u_2$ which originates from a different space.}
}   
\label{proposed_assimilation}
\end{figure}

\textcolor{black}{We illustrate the coarse solution assimilation in Figure (\ref{proposed_assimilation}) and to solve this problem, the information we have can be summarized as follows.}
\begin{enumerate}
    \item A permeability $\Tilde{\kappa}$ which is closed to $\kappa(x)$. If the approximated permeability is given, the correction solution 
    $\Tilde{u}_{H, 2}^n$ ( $\Tilde{\cdot}$ indicates that this is the solution corresponding to $\Tilde{\kappa}(x)$) of this permeability can then be evaluated easily. 
    % \item The multiscale basis associated with $\kappa(x)$; the basis are used to recover the fine solution.
    \item $\{u_{H, 1}^n\}_{n = 1}^M$ where $M<< N$ and $\{\Tilde{u}_{H, 2}^n\}_{n = 1}^N$.
\end{enumerate}

\begin{figure}[H]
\centering
\includegraphics[scale = 0.25]{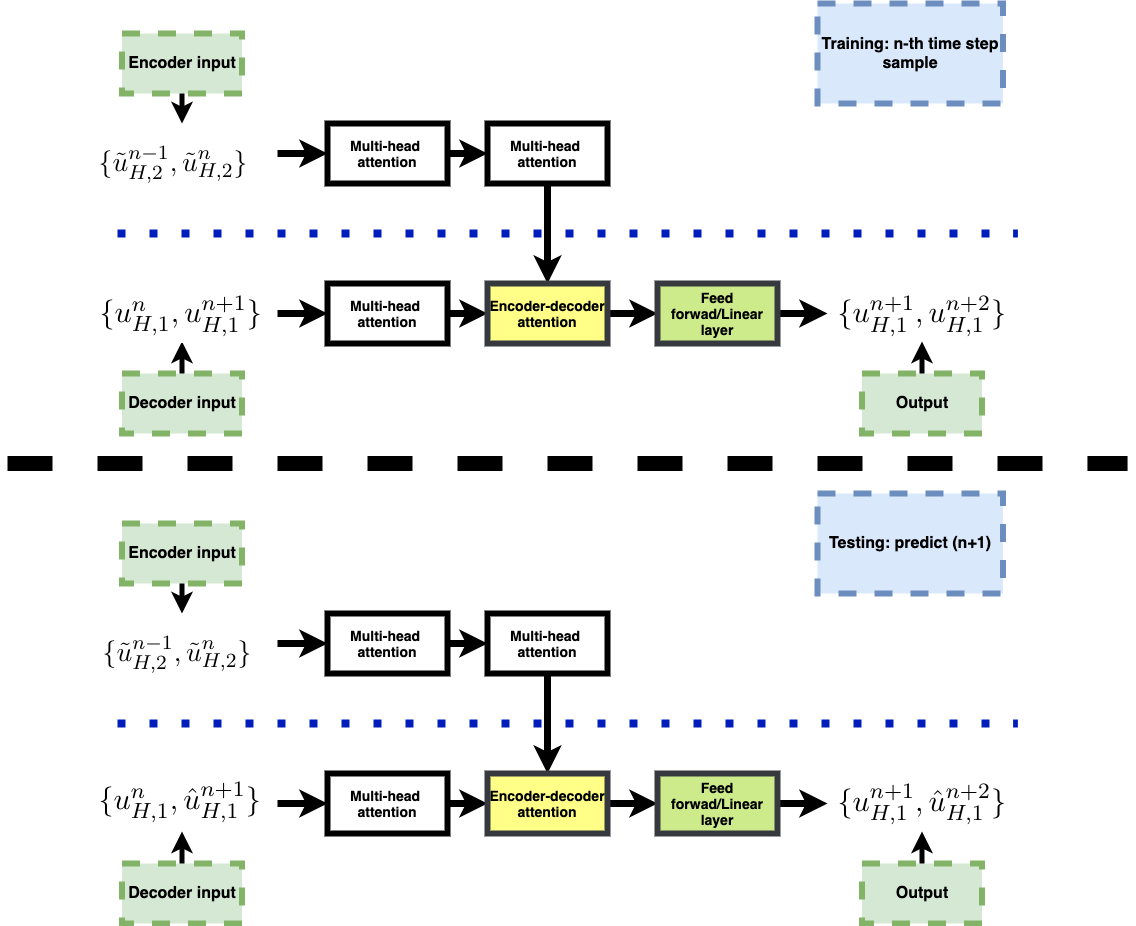}
\caption{Data assimilation model training (top) and testing (bottom).
The inputs of the encoder are $\{\Tilde{u}_{H, 2}^{n-1}, \Tilde{u}_{H, 2}^{n}\}$ which are the solutions of the known equation (with $\Tilde{\kappa}$ as the permeability). The machine model predicts $u_{H, 1}^n$ which can be regarded as a coarse prediction, this solution is then corrected by $\Tilde{u}_{H, 2}^{n}$ which is known and given.
}   
\label{numerical_train_test}
\end{figure}

\textcolor{black}{We have one comment about how to obtain $\{\Tilde{u}_{H, 2}^n\}_{n = 1}^N$. The most accurate approach is to compute it by the splitting scheme. However, this sequence can also be trained and predicted. That is, one can train a model $\mathcal{M}$ with $\{\Tilde{u}_{H, 2}^n\}_{n = 1}^M$ where $M<<N$ and then predicts the rest of the sequence. This process actually is the hybrid explicit-implicit (HEI) learning designed in section (\ref{proposed_sec}) and our numerical experiments show that the prediction is very accurate.}
Another comment is that \textcolor{black}{our proposed hybrid explicit-implicit (HEI) learning method} (training and testing with Transformer) can be used to build this coarse solution assimilation machine learning model.
The training and testing (one example is demonstrated in Figure (\ref{numerical_train_test}) ) are similar to the previous setting except that the input of the encoder will be changed to $\{\Tilde{u}_{H, 2}^n\}_{n = 1}^N$ . 

The advantages of our method (Transformer) are still significant for this application. 
Transformer can handle multiple time series which originate from different spaces.
To be more specific, the model can encode multiple time series and obtain multiple reduced order models of $u_{H, 1}$ which is from the solution space of the target $\kappa$ and $\Tilde{u}_{H, 2}$ which is from the known $\Tilde{\kappa}$.
Moreover, the encoder-decoder attention block is able to learn the connection between two encoded time series, i.e., the model can learn the interplay between two reduced order models. Besides, the model can handle the cumulative error and avoid solution blowing up.

\section{Numerical examples}
\label{numerical_sec}
In this section, we will present some numerical examples.
We solve two equations: linear and nonlinear parabolic equations.
For each equation, we will first present the original example (the model is given) and then discuss the coarse solution assimilation results (the model is unknown).
\textcolor{black}{To generate the training and testing data, we use the splitting scheme (\ref{splitting_splitting_scheme}) for the linear equation and (\ref{splitting_nl_splitting_scheme}) for the nonlinear equation. $\omega$ is set to be $0.5$ for all examples. For the machine learning hyper-parameters, we set $n_e = n_d = 2$, the model hence directly mimics the splitting scheme.}

\subsection{Linear example}
\label{numerical_linear_1_sec}
\begin{figure}[H]
\centering
\includegraphics[scale = 0.45]{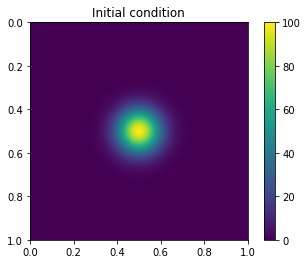}
\caption{Initial condition used in the linear parabolic problems. } 
\label{numerical_exp1_ickappa}
\end{figure}

\begin{figure}[H]
\centering
\includegraphics[scale = 0.45]{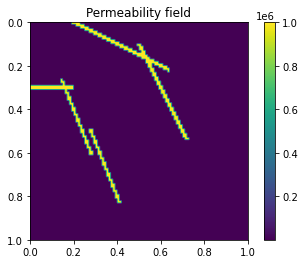}
\includegraphics[scale = 0.45]{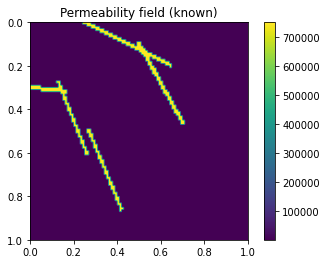}
\caption{ Left: the target permeability. In section (\ref{numerical_linear_1_sec}) (linear) and (\ref{numerical_nl_1_sec}) (nonlinear), this permeability is known. However, we assume it is unknown in the coarse solution assimilation examples. Right: the known permeability $\Tilde{\kappa}$ which is closed to the target unknown permeability $\kappa$ used in the coarse solution assimilation examples. Since $\Tilde{\kappa}(x)$ is given, solution $\{\Tilde{u}_{H, 2}^n\}_{n = 1}^N$ can be derived accordingly by the splitting scheme.} 
\label{numerical_known_kappa}
\end{figure}

In the next two sections, 
we consider the following linear parabolic problem
\begin{align}
    &u_t = \nabla \cdot(\kappa \nabla u ), x\in [0, T]\times \Omega,\\
    & \frac{\partial u}{\partial x} = 0, x\in\partial \Omega,\\
    & u(x, 0) = u_0,
    \label{numerical_eqn1}
\end{align}
where $\Omega = [0, 1]^2$ and $T = 0.00002$, $u_0$ is the initial condition, $\kappa$ is the permeability field, they are demonstrated in Figure (\ref{numerical_exp1_ickappa}) and Figure (\ref{numerical_known_kappa}).

\begin{figure}[H]
\centering
\includegraphics[scale = 0.45]{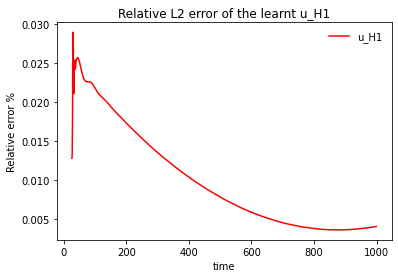}
\includegraphics[scale = 0.45]{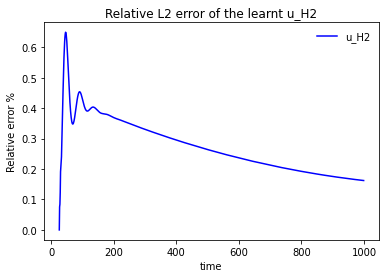}
\caption{Example 1 (linear parabolic). Left, relative error of $u_{H, 1}^n$ with respect to time steps. Right, relative error of $u_{H, 2}^n$ with respect to time steps. The reference solutions are the solutions of the splitting scheme.} 
\label{numerical_exp1_uh12}
\end{figure}

\begin{figure}[H]
\centering
\includegraphics[scale = 0.45]{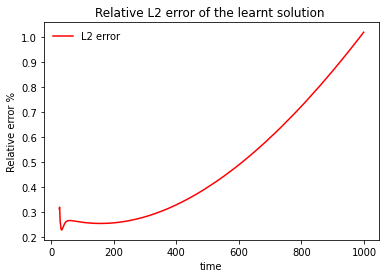}
\includegraphics[scale = 0.45]{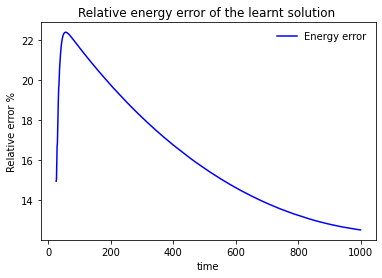}
\caption{Example 1 (linear parabolic) results. Left, the $L_2$ relative error of the solution ($u_{H, 1}^n +u_{H, 2}^n$) with respect to time steps.
Right: the energy relative error of the solution. The reference solution is the FEM solution with mesh $100\times 100$. } 
\label{numerical_exp1_l2}
\end{figure}

The time step $\tau = 0.00000002$ and we train the model with 
the first $25$ steps and predict the remaining $974$ up to the terminal time. We follow the first equation to design the scheme, i.e., $n_d = n_e = 2$. 
To measure the stability of the method, 
we calculate the $L_2$ relative error of $u_{H, 1}^n$ and 
$u_{H, 2}^n$ for all predicted time steps,
the references are the solutions obtained by applying the 
splitting scheme directly. The errors are shown in Figure (\ref{numerical_exp1_uh12}).
We also map the coarse solutions $u_H^n = u_{H, 1}^n +u_{H, 2}^n $ back to the fine scale and compute the $L_2$ relative error, the reference solution is the finite element solution on mesh $100\times 100$ with the backward Euler scheme in time. 
The relative errors are shown in Figure (\ref{numerical_exp1_l2}).

\textcolor{black}{We can see from Figure (\ref{numerical_exp1_uh12}),
the model is able to handle the accumulative error and hence the prediction ($u_{H, 1}^n$) of the model is stable.
Besides, the error of the prediction is also very small (less than $0.03\%$).
$u_{H, 2}^n$ is computed basing on the predictions $u_{H, 1}^n$ of the model,
it is also stable and accurate with relative error less than $0.7\%$.
In Figure (\ref{numerical_exp1_l2}),
we notice that
the relative $L_2$ and energy errors of the solution (the fine mesh solution of the our proposed coarse solution $u_{H}^n$) are also very small and controllable.
Similar results are observed in all following examples and our numerical results are very robust to the hyper-parameters and training, 
we hence conclude that our proposed method is accurate and stable. }

\subsection{Coarse solution assimilation for the linear example}
\label{sec_numerical_assimilation}
\begin{figure}[H]
\centering
\includegraphics[scale = 0.45]{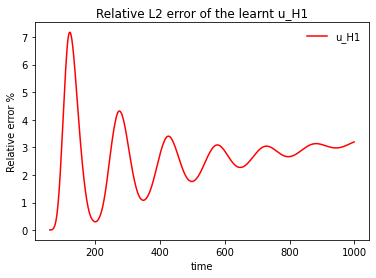}
\caption{ Results of the linear coarse solution assimilation example. Relative error of $u_{H, 1}^n$ with respect to time steps.
} 
\label{numerical_assimilation_uh1}
\end{figure}

To demonstrate our idea, we will solve the same linear equation (\ref{numerical_eqn2}). 
We will train the model with the first $60$ steps and predicts the remaining $939$ steps.
All the other settings (eg, source, time step...) are the same as before.
We calculate $L_2$ relative error of $u_{H, 1}^n$ for all predicted time steps.
Please note that there is no need to compute $u_{H, 2}^n$ since the $u_{H, 1}^n$ is corrected by the given $\Tilde{u}_{H, 2}^n$ instead.
The result is presented in Figure (\ref{numerical_assimilation_uh1}).

\begin{figure}[H]
\centering
\includegraphics[scale = 0.45]{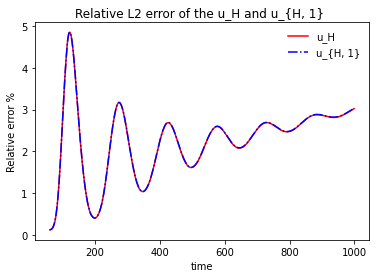}
\includegraphics[scale = 0.45]{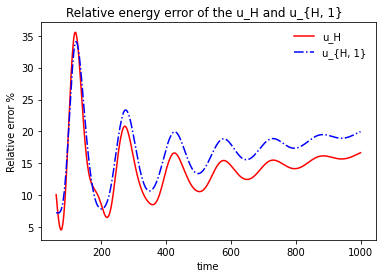}
\caption{ Results of the linear coarse solution assimilation example. Left: relative $L_2$ error with respect to time steps. Right: 
relative energy error with respect to time steps. Red color dashed line is the fine scale $u_{H}^n$ (which maps the coarse solution to the fine scale mesh), while the blue curve is the fine scale $u_{H, 1}^n$.
} 
\label{numerical_assimilation_errors}
\end{figure}

We also map the coarse solutions $u_{H, 1}^N$ and $u_H^N = u_{H, 1}^N +\Tilde{u}_{H, 2}^N $ back to the fine scale and compute the $L_2$ and energy relative errors. The idea is to correct $u_{H, 1}^N$ with $\Tilde{u}_{H, 2}^N$ which is the correction of the approximated permeability.
The results are shown in Figure (\ref{numerical_assimilation_errors}) and please check Table (\ref{numerical_assimilation_table}) for the details.

\begin{table}[H]
\centering
\begin{tabular}{||c c c||} 
\hline
Coarse & $L_2$  & Energy \\ [0.5ex] 
\hline
$u_{H, 1}$  & 2.262591\% & 17.142532\% \\ [0.5ex]
\hline
$u_{H}$  & 2.25732\% & 14.435881\%  \\ [0.5ex]
\hline
\end{tabular}
\caption{Linear coarse solution assimilation example relative errors. We here present the $L_2$ and the energy relative errors of the solutions obtained by mapping $u_{H, 1}^n$ and $u_{H}^n$ to the fine scale mesh. The errors are calculated as the average over all time steps.}
\label{numerical_assimilation_table}
\end{table}

\textcolor{black}{
We can observe from Figure (\ref{numerical_assimilation_uh1}) and Figure (\ref{numerical_assimilation_errors})
both prediction and the solution are stable and are accurate.
Since the target $\kappa$ is not given,
the model is supposed to learn $u_{H, 1}$ basing on $\Tilde{u}_{H, 2}$, moreover, $\Tilde{u}_{H, 2}$ should correct and improve $u_{H, 1}$. This can be seen from Table (\ref{numerical_assimilation_table}). 
The  $L_2$ and energy errors are small
and adding $\Tilde{u}_{H, 2}$ indeed improves the relative energy error of the solution.
}

\textcolor{black}{We want to give a comment before finishing the linear coarse solution assimilation example. 
It may be costly to pre-compute all $\Tilde{u}_{H, 2}^n$. We provide two strategies which can reduce the computational cost. 
Firstly, one can calculate $\Tilde{u}_{H, 2}^n$ by following the scheme (\ref{splitting_assimilation_2}).
Compared to pre-compute all $\Tilde{u}_{H, 2}^n$, one avoids the computation of $\Tilde{u}_{H, 1}^n$ which is the contrast dependent solution of the approximation model $\Tilde{\kappa}$.
Secondly, we can train a model that predicts $\Tilde{u}_{H, 1}^n$ and calculate $\Tilde{u}_{H, 2}^n$ of the approximated model $\Tilde{\kappa}$, this is exactly the hybrid explicit-implicit (HEI) learning in the section (\ref{numerical_linear_1_sec}). Since the maximal relative error of $u_{H, 2}^n$ is less than $0.6\%$ in the last example (\ref{numerical_linear_1_sec}), we believe our method can provide an accurate $\Tilde{u}_{H, 2}^n$ for later use in the assimilation. In all, this strategy consists of two stages. In the first stage, one can train a model for $\Tilde{\kappa}$ and then predict and compute all $\Tilde{u}_{H, 1}^n$ and $\Tilde{u}_{H, 2}^n$. In the second stage, one can train a different model which relies on the pre-computed $\Tilde{u}_{H, 2}^n$ to perform the coarse data assimilation.
}

\subsection{Nonlinear example}
\label{numerical_nl_1_sec}
\begin{figure}[H]
\centering
\includegraphics[scale = 0.45]{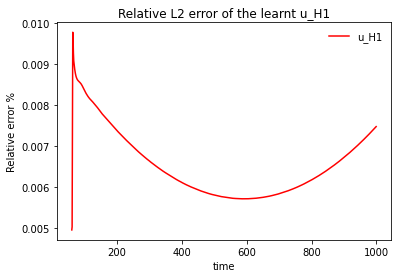}
\includegraphics[scale = 0.45]{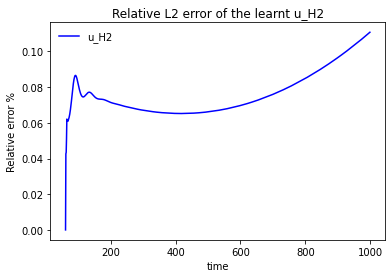}
\caption{ Results of the nonlinear example. Left: the $L_2$ relative error of $u_{H, 1}^n$ with respect to time steps. Right: 
the $L_2$ relative error of the solution $u_{H, 2}^n$ with respect to time steps.
} 
\label{numerical_nonlinear_uh12}
\end{figure}

\begin{figure}[H]
\centering
\includegraphics[scale = 0.45]{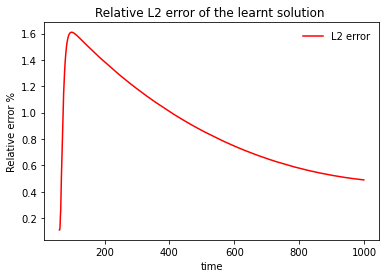}
\includegraphics[scale = 0.45]{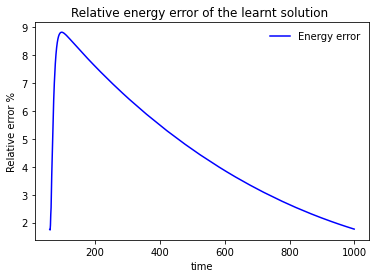}
\caption{Example 3 (nonlinear parabolic problem) results. Left, the $L_2$ relative error of the solution ($u_{H, 1}^n +u_{H, 2}^n$) with respect to time steps.
Right: the energy relative error of the solution. The reference solution is the FEM solution with mesh $100\times 100$. } 
\label{numerical_exp3_l2}
\end{figure}

In the next two sections, 
we consider the following nonlinear parabolic problem:
\begin{align}
    &u_t = \nabla \cdot(e^{u}\kappa_0(x) \nabla u ), x\in [0, T]\times \Omega,\\
    & \frac{\partial u}{\partial x} = 0, x\in\partial \Omega,\\
    & u(x, 0) = u_0,
    \label{numerical_eqn2}
\end{align}
where $\Omega = [0, 1]^2$ and $T = 0.00002$, $u_0 = e^{(-(x_1-0.5)^2-(x_2-0.5)^2)/0.01}$ is the initial condition, $\kappa_0$ is the permeability field, they are demonstrated in Figure (\ref{numerical_exp1_ickappa}) and Figure (\ref{numerical_known_kappa}).
We set $\tau = 0.00000002$ and we train the model with the first 60 time steps and predict the remaining 939 steps.
We calculate the $L_2$ relative error of $u_{H, 1}^n$ and $u_{H, 2}^n$ for all predicted time steps, the references are the solutions obtained by applying the splitting scheme directly.  The errors are shown in Figure (\ref{numerical_nonlinear_uh12}).
We also map the coarse solutions $u_H^n = u_{H, 1}^n +u_{H, 2}^n $ back to the fine scale and compute the $L_2$ and energy relative error, the reference solution is the finite element solution on mesh $100\times 100$ with the backward Euler scheme in time. 
The relative errors are shown in Figure (\ref{numerical_exp3_l2}).

\textcolor{black}{Similar with the linear example, the nonlinear model is also
accurate and stable (please check Figure (\ref{numerical_nonlinear_uh12}) left):
the relative error in the prediction of $u_{H, 1}$ is less than 
$0.01\%$, the computed $u_{H, 2}$ which relies on the prediction of the model is also small (less than $0.15\%$).
If we look at Figure (\ref{numerical_exp3_l2}), the fine scale solution which is obtained by mapping $u_{H, 1}^n+u_{H, 2}^n$ back to the fine scale mesh is also stable, the relative errors also show that our method is very accurate.}

\subsection{Coarse solution assimilation for the nonlinear example}
\label{sec_numerical_nl_assimilation}
\begin{figure}[H]
\centering
\includegraphics[scale = 0.45]{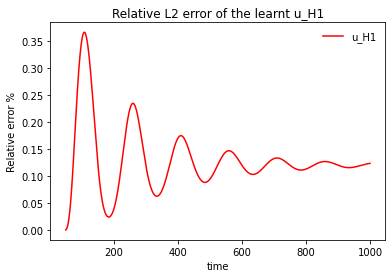}
\caption{ Results of the nonlinear coarse solution assimilation example. The $L_2$ relative error of $u_{H, 1}^n$ with respect to time steps.
} 
\label{numerical_nl_assimilation_uh1}
\end{figure}

To demonstrate our idea, we will solve the same nonlinear equation (\ref{numerical_eqn2}). 
We will train the model with the first $50$ steps and predicts the remaining $949$ steps,
all the other settings (eg, source, time step...) are the same as before.
We calculate $L_2$ relative error of $u_{H, 1}^n$ for all predicted time steps.
Please note that there is no need to compute $u_{H, 2}^n$ since the $u_{H, 1}^n$ is corrected by the given $\Tilde{u}_{H, 2}^n$.
The result is presented in Figure (\ref{numerical_nl_assimilation_uh1}).

\begin{figure}[H]
\centering
\includegraphics[scale = 0.45]{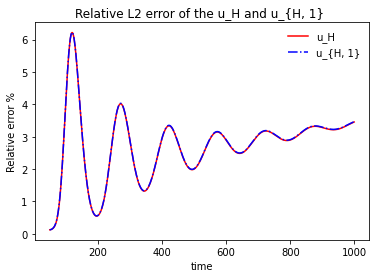}
\includegraphics[scale = 0.45]{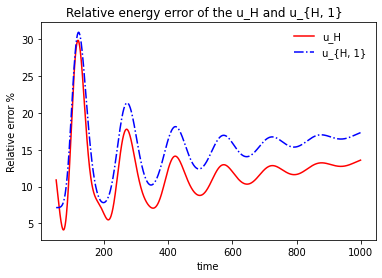}
\caption{ Results of the nonlinear coarse solution assimilation example. Left: relative $L_2$ error with respect to time steps. Right: 
relative energy error with respect to time steps. Red color dashed line is the fine scale $u_{H}$, while the blue curve is the fine scale $u_{H,1}$.
} 
\label{numerical_nl_assimilation_errors}
\end{figure}

We also map the coarse solutions $u_{H, 1}^n$ and $u_H^n = u_{H, 1}^n +\Tilde{u}_{H, 2}^n $ back to the fine scale and compute the $L_2$ and energy relative errors. The idea is to correct $u_{H, 1}^n$ with $\Tilde{u}_{H, 2}^n$ which is the correction of the approximated permeability.
The results are shown in Figure (\ref{numerical_nl_assimilation_errors}) and please check Table (\ref{numerical_nl_assimilation_table}) for the details.

\begin{table}[H]
\centering
\begin{tabular}{||c c c||} 
\hline
Coarse & $L_2$  & Energy \\ [0.5ex] 
\hline
$u_{H, 1}$  & 2.742292\% & 15.498801\% \\ [0.5ex]
\hline
$u_{H}$  & 2.731022\% & 11.996128\%  \\ [0.5ex]
\hline
\end{tabular}
\caption{Nonlinear coarse solution assimilation example relative errors. We here present the $L_2$ and the energy relative errors of the solutions obtained by mapping $u_{H, 1}$ and $u_{H}$ to the fine scale mesh.  The errors are calculated as the average over all time steps.}
\label{numerical_nl_assimilation_table}
\end{table}

\textcolor{black}{
Figure (\ref{numerical_nl_assimilation_uh1}) demonstrates the 
error of the model prediction. The relative error is less than $0.4\%$ and is very stable.
The  errors are shown in Figure (\ref{numerical_nl_assimilation_errors}), we can see that the solution is very stable and accurate. 
If we look at Table (\ref{numerical_nl_assimilation_table}),
the $\Tilde{u}_{H, 2}$ indeed corrects $u_{H, 1}$ and improves the solution. To reduce the cost of precomputing $\Tilde{u}_{H, 2}^n$, one can use the same hybrid explicit-implicit (HEI) learning provided in section (\ref{sec_numerical_assimilation}). In particular for the strategy two which trains a separate model for $\Tilde{u}_{H, 2}^n$, the maximal relative error of $u_{H, 2}^n$ in example \ref{numerical_nl_1_sec} is less than $0.1\%$, we believe that the trained model can provide an accurate $u_{H, 2}^n$ for the future use in the coarse solution assimilation.}

\section{Conclusion}
In this work, we design a \textcolor{black}{hybrid explicit-implicit (HEI) learning} approach to accelerate the contrast-independent splitting scheme and solve coarse data assimilation problems.
We proposed three strategies, i.e., one can either learn 
$u_1$ which represents the main feature of the problem, or $u_2$
which is the contrast independent part of the solution, or both $u_1$ and $u_2$.
We focus on learning $u_1$ because the computation of $u_1$ is implicit and time consuming if we follow the standard splitting scheme.
We interpret the splitting scheme as a language translation problem and hence use the Transformer to solve it.
Transformer is powerful in our problem since it can encode the inputs of different spaces and reduce the accumulative error.
We also find that the splitting is a great platform to solve the coarse solution assimilation problem. All three strategies can be used and our machine learning model is also powerful in solving the coarse solution assimilation problems. We demonstrate and verify our proposed method with four examples, and the results show that our method is stable and accurate.
In the future, we will study the other two strategies and improve our model.

\section{Acknowledgement}
Guang Lin gratefully acknowledges the support of the National Science Foundation (DMS-1555072, DMS-1736364, DMS-2053746, and DMS-2134209), and Brookhaven National Laboratory Subcontract 382247, and U.S. Department of Energy (DOE) Office of Science Advanced Scientific Computing Research program DE-SC0021142).

\bibliographystyle{abbrv}
\bibliography{references, references1, references2, references3, references4}
\end{document}